\documentclass[12pt]{article}
\usepackage{setspace,graphicx}
\usepackage{amsmath,amsthm,amsfonts,amssymb}
\usepackage{xcolor}
\usepackage[colorlinks,citecolor=blue,urlcolor=blue,linkcolor=blue]{hyperref}
\hypersetup{pdfauthor={Alexey Kuznetsov, Minjian Yuan},
  pdftitle={Rational extensions of the Ornstein--Uhlenbeck process},
  pdfkeywords={Ornstein--Uhlenbeck process, exceptional Hermite polynomials, Mehler formula, Darboux transformation, transition density, quantum harmonic oscillator},
  pdfpagemode=UseNone}
\def\r{\mathbb R}
\def\i{{\mathrm i}}
\newcommand{\z}{\mathbb Z}
\newcommand{\n}{\mathbb N}
\def\d{\mathrm d}
\def\k{\textbf{\textsl{k}}}
\def\j{\textbf{\textsl{j}}}
\DeclareMathOperator{\Wr}{Wr}
\newtheorem{theorem}{Theorem}
\newtheorem{lemma}{Lemma}
\newtheorem{proposition}{Proposition}

\newtheorem{definition}{Definition}

\title{Rational extensions of the Ornstein--Uhlenbeck process}
\author{Alexey Kuznetsov, Minjian Yuan\footnote{Dept. of Mathematics and Statistics,
York University, 4700 Keele Street, Toronto, ON, M3J 1P3, Canada.
\\ Email: akuznets@yorku.ca, yuanm@yorku.ca}}
\date{\today}

\begin{document}
\maketitle
\begin{abstract}
We introduce and classify rational extensions of the Ornstein--Uhlenbeck process that are solvable by polynomials. These diffusions are closely related to exceptional Hermite polynomials and rational extensions of the quantum harmonic oscillator. We obtain explicit formulas for their transition densities, fundamental solutions and Green functions. We also show that each such diffusion can be constructed from killed Brownian motion by a finite sequence of Darboux transformations followed by a Doob $h$-transform. Finally, we derive an explicit formula for the connection polynomials in the generalized Mehler formula and use it to prove several conjectures concerning these polynomials.
\end{abstract}



\section{Introduction and main results}

The Ornstein--Uhlenbeck (OU) process, studied by Uhlenbeck and
Ornstein in 1930 \cite{UhlenbeckOrnstein1930}, originated in the
theory of Brownian motion. It models the velocity of a particle
subject to linear friction and random forcing, following
Langevin's 1908 equation of motion \cite{Langevin1908}.
In the normalization used here, the process
$Y=\{Y_t\}_{t\geqslant 0}$ satisfies the stochastic differential
equation (SDE)
\[
\d Y_t = -Y_t\,\d t + \d W_t,
\qquad Y_0=x,
\]
where $W$ is a standard Brownian motion.
The OU process is particularly useful because of its mathematical
tractability. Its transition probability density with respect to
Lebesgue measure is
\begin{equation}\label{p_t^OU_formula}
p_t^{\mathrm{OU}}(x,y)
=
\frac{1}{\sqrt{\pi(1-e^{-2t})}}
\exp\left(
-\frac{(y-xe^{-t})^2}{1-e^{-2t}}
\right),
\qquad t>0,\quad x,y\in\r.
\end{equation}
Its invariant distribution is Gaussian with  mean zero and variance  $1/2$, and the Hermite
polynomials $H_n$ form an orthogonal basis of eigenfunctions
of its infinitesimal generator in the $L_2$ space associated
with this invariant measure
\cite{Bogachev2018,Borodin2002,Pavliotis2014}.

OU processes and their extensions have applications across
physics, biology, climate science, and finance.
In turbulence, OU velocity models and their extensions are
used in Lagrangian descriptions of particle dispersion
\cite{Thomson1987}.
In active-matter models, OU processes represent exponentially
correlated self-propulsion velocities \cite{Fodor2016}.
In chemical kinetics and systems biology, multivariate OU
processes describe fluctuations about stable deterministic
steady states of reaction networks within the linear-noise
approximation \cite{ElfEhrenberg2003,VanKampen2007}.
In climate science, OU models describe slow climate variables,
such as sea-surface temperature anomalies, driven by rapid
weather fluctuations and subject to linear stabilizing feedback
\cite{Hasselmann1976,FrankignoulHasselmann1977}. Further biological applications include phylogenetic models of
quantitative-trait evolution toward adaptive optima under
stabilizing selection \cite{Hansen1997}, as well as models of
animal movement and home-range estimation based on OU processes
and their extensions \cite{Fleming2014}.
In mathematical neuroscience, OU processes model subthreshold
membrane potentials, with neuronal firing described by
first-passage times through a threshold
\cite{RicciardiSacerdote1979}.
In mathematical finance, the OU process models the instantaneous
short rate in Vasicek's term-structure model \cite{Vasicek1977}.

The OU process is closely connected to the quantum harmonic
oscillator. Its infinitesimal generator
$\mathcal L^{\mathrm{OU}}$
is related to the harmonic oscillator Hamiltonian
$\mathcal H:=-\partial_x^2+x^2$ (see \cite[Section 8.3]{Teschl2014})
by the ground-state transformation
\cite[Section 4.4]{Pavliotis2014}
\[
\mathcal L^{\mathrm{OU}}:=\frac12\partial_x^2-x\partial_x=-\frac12 e^{x^2/2}(\mathcal H-1)e^{-x^2/2}.
\]
This correspondence relates the Hermite polynomials
$H_n(x):=(-1)^n e^{x^2}\partial_x^n e^{-x^2}$
to the Hermite functions $\psi_n(x):=e^{-x^2/2}H_n(x)$:
they satisfy $\mathcal L^{\mathrm{OU}}H_n=-nH_n$ and
$\mathcal H\psi_n=(2n+1)\psi_n$, respectively.
It also provides a natural framework for constructing rational
extensions of the OU process. Rational extensions of the harmonic
oscillator that are exactly solvable by polynomials have been
classified via Darboux transformations \cite{Gomez-Ullate_2014},
and explicit formulas for their propagators are available
\cite{Pupasov_2015}.
One of our goals is to use these results to construct and study
rational extensions of the OU process, which we define next. 

Let $R=P/Q$ be a rational function without real poles, where
$P,Q\in\r[x]$ and $\deg P<\deg Q$.
A {\it rational extension} of the OU process is a diffusion $X$ on $\r$
satisfying
\[
\d X_t=-(X_t+R(X_t))\,\d t+\d W_t.
\]
Since $R$ and $R'$ are bounded, the drift is globally Lipschitz
and has at most linear growth. Thus, for every initial value
$X_0=x\in\r$, this SDE has a unique global strong solution
\cite[Theorem~5.2.1]{Oksendal2003}.

The infinitesimal generator of $X$ is
\begin{equation}\label{def_L_R}
\mathcal L_R:=\frac12\partial_x^2-(x+R(x))\partial_x.
\end{equation}
Define the weight
\[
\mu_R(x):=\exp\left(-x^2-2\int_0^xR(y)\,\d y\right).
\]
Since $R(x)=O(|x|^{-1})$, we have
$\mu_R(x)=O(|x|^\alpha e^{-x^2})$ as $|x|\to\infty$
for some $\alpha>0$, and hence $\mu_R\in L_1(\r,\d x)$.
Conjugation by $\sqrt{\mu_R}$ maps $-2\mathcal L_R$
unitarily to a Schrödinger operator
\[
\sqrt{\mu_R}\, (-2 \mathcal L_R)\, \frac{1}{\sqrt{\mu_R}} = - \partial_x^2+ (x+R(x))^2-1-R'(x).
\]
The potential of this Schrödinger operator is $x^2+O(1)$ as $|x|\to \infty$. The Kato--Rellich theorem \cite[Theorem~6.4]{Teschl2014} implies $\mathcal L_R$ is essentially self-adjoint
on $C_c^\infty(\r)$ in $L_2(\r,\mu_R(x)\d x)$.
Its closure, denoted by the same symbol, has purely discrete
spectrum and an orthonormal basis of eigenfunctions
\cite[Theorem~XIII.67]{ReedSimon1978}. These eigenfunctions play an important role in the following definition, which we adapted from
\cite[Definition~1.1]{Gomez-Ullate_2014}.
\begin{definition}
A rational extension of the OU process is \emph{solvable by
polynomials}, or simply \emph{solvable}, if there exists a fixed
function $w\not\equiv0$ such that, for all but finitely many integers
$n \geqslant 0$, there is a polynomial $p_n\in\r[x]$ of degree $n$
for which $\phi_n=wp_n$ is an eigenfunction of $\mathcal L_R$
in $L_2(\r,\mu_R(x)\d x)$.
\end{definition}

We introduce notation for our classification of solvable rational
extensions of the OU process. We write $\n=\{1,2,\ldots\}$ and
$\z_{\geqslant0}=\{0\} \cup \n$.
The Wronskian of smooth functions $f_1,\ldots,f_m$ is
\[
\Wr[f_1,\ldots,f_m]:=\det[\partial_x^{i-1}f_j]_{1\leqslant i,j\leqslant m}.
\]
For $\j=(j_1,\ldots,j_m)\in\z_{\geqslant0}^m$, write
$H_{\j}(x):=\Wr[H_{j_1},\ldots,H_{j_m}](x)$.
Throughout, $\k=(k_1,\ldots,k_\ell)\in\n^\ell$ denotes
a strictly increasing sequence of length $\ell$, and we write
$\k-1:=(k_1-1,\ldots,k_\ell-1)$. We use the convention $H_\varnothing=1$ and interpret empty products as $1$.

Such a sequence $\k$ is called a \emph{Krein--Adler sequence}
if $\ell\geqslant2$ is even and
$k_{2i}=k_{2i-1}+1$ for $1\leqslant i\leqslant\ell/2$.
For every Krein--Adler sequence $\k$, the polynomials
$H_{\k}$ and $H_{\k-1}$ are even and strictly positive on $\r$
(see Lemma~\ref{Lemma_H_k_H_k-1}).
The following theorem is a diffusion analogue of
\cite[Theorem~1.1]{Gomez-Ullate_2014}, which classifies
solvable regular rational extensions of the harmonic oscillator.

\begin{theorem}\label{thm:classification}
A rational extension of the OU process is solvable if and only
if either $R\equiv0$ or there exists a Krein--Adler sequence $\k$
such that
\begin{equation}\label{def_R_k}
R(x)=\frac{H_{\k}'(x)}{H_{\k}(x)}-\frac{H_{\k-1}'(x)}{H_{\k-1}(x)}.
\end{equation}
\end{theorem}

For a Krein--Adler sequence $\k$, let $X^{\k}$ denote the
corresponding diffusion, with infinitesimal generator
\[
\mathcal L^{\k}:=\frac12\partial_x^2-\left(x+\frac{H_{\k}'(x)}{H_{\k}(x)}-\frac{H_{\k-1}'(x)}{H_{\k-1}(x)}\right)\partial_x.
\]
For $n\in\z_{\geqslant0}$, write
$H_{\k,n}(x):=\Wr[H_{k_1},\ldots,H_{k_\ell},H_n](x)$.
The nonzero members of this family, corresponding to $n\notin\k$,
are called \emph{exceptional Hermite polynomials}.
Their differential equations, orthogonality, and completeness
yield the following result
\cite{Duran_2014,Gomez-Ullate_2014,Gomez-Ullate_2020}.

\begin{theorem}\label{thm:spectrum}
Let $\k$ be a Krein--Adler sequence.
The rational functions
$\varphi_{\k,n}:=H_{\k,n}/H_{\k-1}$,
$n\in\z_{\geqslant0}\setminus\k$, satisfy
\[
\mathcal L^{\k}\varphi_{\k,n}=-n\varphi_{\k,n}.
\]
They form an orthogonal basis of $L_2(\r,\nu_{\k}(x)\d x)$, where
\[
\nu_{\k}(x):=e^{-x^2}\frac{H_{\k-1}(x)^2}{H_{\k}(x)^2}.
\]
For $n,m\in\z_{\geqslant0}\setminus\k$,
\begin{equation}\label{eq:norm}
\int_{\r}\varphi_{\k,n}(x)\varphi_{\k,m}(x)\nu_{\k}(x)\,\d x=\delta_{n,m}\sqrt\pi\,2^{\ell+n}n!\prod_{j=1}^{\ell}(n-k_j).
\end{equation}
\end{theorem}

Our next goal is to determine the transition probability density of $X^{\k}$. The main tool is the generalized Mehler formula, which states that for every strictly increasing sequence $\k \in {\mathbb N}^{\ell}$ there exist {\emph{connection polynomials}} $Q_n^{\k}(x,y)$, $0\leqslant n \leqslant k_{\ell}+1$, such that for all $z,x,y\in\mathbb C$ with $|z|<1$,
\begin{equation}\label{eqn:Mehler}
\sum_{n\in\z_{\geqslant0}\setminus\k}\frac{z^n H_{\k,n}(x)H_{\k,n}(y)}{2^{n+\ell}n!\prod_{j=1}^{\ell}(n-k_j)}
=\frac{1}{\sqrt{1-z^2}}\exp\!\left(\frac{2zxy-z^2(x^2+y^2)}{1-z^2}\right)\sum_{n=0}^{k_\ell+1}z^n Q_n^{\k}(x,y).
\end{equation} 
This formula was first established for Krein--Adler sequences
by Pupasov-Maksimov \cite{Pupasov_2015} and extended to arbitrary
strictly increasing sequences in \cite{Kuznetsov2026}. By \cite[Proposition 1.3]{Kuznetsov2026}, we have
\begin{equation}\label{eqn_Qk0}
Q_0^{\k}(x,y)=(-2)^\ell\Big[\prod_{j=1}^{\ell}k_j\Big] H_{\k-1}(x)H_{\k-1}(y).
\end{equation}
The remaining polynomials $Q_n^\k(x,y)$,
$1\leqslant n\leqslant k_\ell+1$, can be computed recursively using \cite[Eq.~(1.7)]{Kuznetsov2026}.
Our next result gives a closed-form expression for the connection polynomials.

\begin{theorem}\label{thm:Mehler_kernel}
Let $\k\in\n^\ell$ be strictly increasing, with $\ell\geqslant 1$.
Set $\k\setminus k_i:=(k_j)_{1\leqslant j\leqslant \ell,\ j\ne i}$ and  $T_n(x,y):=H_n(x)H_n(y)/(2^n n!)$ for $n\geqslant 0$
and $T_n(x,y):=0$ for $n<0$.
Then, for $1\leqslant n\leqslant k_\ell+1$,
\begin{equation}\label{eq:explicit}
Q_n^{\k}(x,y)=2^{\ell-1}\sum_{i=1}^{\ell}2^{k_i}k_i!\Big[\prod_{j\ne i}(k_i-k_j)\Big]H_{\k\setminus k_i}(x)H_{\k\setminus k_i}(y)\big(T_{k_i+1-n}(x,y)-T_{k_i-1-n}(x,y)\big).
\end{equation}
\end{theorem}

In Section~\ref{Section_Q_n}, we use the explicit formula
\eqref{eq:explicit} to prove Conjectures 1.7 and 1.8 and the first assertion of Conjecture 1.9 in
\cite{Kuznetsov2026}. In the next theorem we use the connection polynomials to obtain explicitly
the transition density of $X^{\k}$. 

\begin{theorem}\label{thm:transition-density}
Let $\k$ be a Krein--Adler sequence. For $t>0$ and $x,y\in\r$,
the transition density of $X^{\k}$ with respect to Lebesgue
measure is
\begin{equation}\label{eqn:p_t^X}
p_t^{\k}(x,y)=p_t^{\mathrm{OU}}(x,y)\frac{H_{\k-1}(y)}{H_{\k}(y)^2H_{\k-1}(x)}\sum_{n=0}^{k_\ell+1}e^{-nt}Q_n^{\k}(x,y).
\end{equation}
\end{theorem}

The process $X^{\k}$ is positive recurrent, with invariant
probability density
\begin{equation}\label{def_pi_k}
\pi_{\k}(x)=\frac{2^\ell}{\sqrt{\pi}}
\Big[\prod_{j=1}^{\ell}k_j\Big]\, e^{-x^2}\frac{H_{\k-1}(x)^2}{H_{\k}(x)^2}.
\end{equation}
The normalizing constant follows from \eqref{eq:norm} with
$n=m=0$. Since $R$ in \eqref{def_R_k} is odd
(see Lemma~\ref{Lemma_H_k_H_k-1}), the law
of $X^{\k}$ under ${\mathbb P}_x$ coincides with the law of $-X^{\k}$ under ${\mathbb P}_{-x}$.

We next determine the fundamental solutions and the Green
function of $X^{\k}$. For $\lambda>0$, the fundamental solutions
are positive increasing and decreasing solutions of
$\mathcal L^{\k}f=\lambda f$, determined up to positive multiplicative  constants. For the OU process, the fundamental solutions are
$g_{\lambda}^{\pm}(x)=F_{\lambda}(\pm x)$, where
\[
F_{\lambda}(x):=e^{x^2/2}D_{-\lambda}(-\sqrt{2}x),
\]
and $D_p$ is the parabolic cylinder function
\cite[Section~12.2]{NIST}. The identities
\cite[Eqs.~(12.2.5) and~(12.2.11)]{NIST} give
\begin{equation}\label{g^pm_lambda_Wronskian}
\Wr[g_{\lambda}^-,g_{\lambda}^+](x)
=\frac{2\sqrt{\pi}}{\Gamma(\lambda)}e^{x^2}.
\end{equation}
The following theorem gives the corresponding fundamental solutions
for $X^{\k}$.

\begin{theorem}\label{thm:rational-ou-fundamental}
Let $\k$ be a Krein--Adler sequence and $\lambda>0$. Define
\begin{equation}\label{def_F_k_lambda}
F_{\k,\lambda}(x)
:=\Wr[H_{k_1},\ldots,H_{k_\ell},F_\lambda](x).
\end{equation}
The increasing and decreasing fundamental solutions of
$\mathcal L^{\k}f=\lambda f$ can be chosen as
\begin{equation}\label{def_f_lambda_pm}
f_\lambda^+(x):=\frac{F_{\k,\lambda}(x)}{H_{\k-1}(x)},
\qquad f_\lambda^-(x):=f_\lambda^+(-x).
\end{equation}
Their Wronskian is
\begin{equation}\label{f^pm_lambda_Wronskian}
\Wr[f_\lambda^-,f_\lambda^+](x)
=\frac{\sqrt{\pi}\,2^{\ell+1}\prod_{j=1}^{\ell}(\lambda+k_j)}
{\Gamma(\lambda)\nu_{\k}(x)}.
\end{equation}
\end{theorem}

Fundamental solutions play a central role in the study of diffusion
processes
(see \cite{Borodin2002}). We give two applications. Define the running maximum $M_t:=\max\{ X^{\k}_s \; : \; 0\leqslant s \leqslant t\}$ and the first-passage time $h_y:=\inf\{ t\geqslant 0 \; : \; X^{\k}_t=y\}$. 
Then, if $y\geqslant x$ and $T$ is independent of $X^{\k}$ and has an exponential distribution with rate $\lambda>0$, we have 
\[
{\mathbb P}_x(M_T \geqslant y)={\mathbb E}_x \big[e^{-\lambda h_y}\big]=\frac{f^+_{\lambda}(x)}{f^+_{\lambda}(y)}.
\]
There are similar formulas involving the running minimum of $X^{\k}$ and $f^-_{\lambda}$ in the case $y\leqslant x$. See \S II.10 and \S II.19 in \cite{Borodin2002}. 
Next, Theorem~\ref{thm:rational-ou-fundamental} and the standard
Green-function formula (see 
\cite[page 19]{Borodin2002}) yield, for
$\lambda>0$ and $x,y\in\r$,
\[
G_\lambda^{\k}(x,y)
:=\int_0^\infty
e^{-\lambda t}p_t^{\k}(x,y)\,\d t=\frac{\Gamma(\lambda)\nu_{\k}(y)}
{\sqrt{\pi}\,2^\ell\prod_{j=1}^{\ell}(\lambda+k_j)}
f_\lambda^+(x\wedge y)f_\lambda^-(x\vee y),
\]
where $x\wedge y=\min\{x,y\}$ and $x\vee y=\max\{x,y\}$.

We next give a probabilistic construction of the solvable rational
extensions classified in Theorem~\ref{thm:classification}.
The natural starting point is Brownian motion killed at rate
$x^2/2$: its generator is $-\mathcal H/2$, and its Doob
$h$-transform with $h(x)=e^{-x^2/2}$ is the OU process. In \cite{KuznetsovYuan}, we introduced \emph{Darboux transformations
of diffusion processes}. Under suitable conditions, this construction
takes a Brownian motion $Z$ killed at rate $c(x)$ and a positive
$\lambda$-invariant function $h$, and produces a Brownian motion
$\widetilde Z$ killed at a new rate $\tilde c(x)$. The transition
density of $\widetilde Z$ is expressed explicitly in terms of
that of $Z$ \cite[Theorem~4.8]{KuznetsovYuan}. The following
result shows that iterating this construction, followed by a
Doob $h$-transform, yields every nontrivial solvable rational
extension of the OU process.

\begin{theorem}\label{thm:OU-Darb}
Let $\k$ be a Krein--Adler sequence. The process $X^{\k}$ can be
obtained from Brownian motion killed at rate $x^2/2$ by
$k_{\ell}+1-\ell$ successive Darboux transformations of diffusion
processes in the sense of \cite{KuznetsovYuan}, followed by a
Doob $h$-transform.
\end{theorem}

Now we present several probabilistic applications of our results.

\begin{proposition}\label{cor:OU-bridge}
Let $\k$ be a Krein--Adler sequence.
For $t>0$ and $x,y\in\r$, let
$\mathbb E_{x\to y}^{\mathrm{OU},t}$ denote expectation
under the OU bridge from $x$ to $y$ over $[0,t]$.
Then
\begin{equation}\label{eq:OU-bridge-functional}
\mathbb E_{x\to y}^{\mathrm{OU},t}
\left[
\exp\left(
\int_0^t(\log H_\k)''(Y_s)\,\d s
\right)
\right]
=
\frac{1}{H_\k(x)H_\k(y)}
\sum_{n=0}^{k_\ell+1}e^{(\ell-n)t}Q_n^\k(x,y).
\end{equation}
\end{proposition}

For $\k=(1,2)$, this identity becomes
\[
\mathbb E_{x\to y}^{\mathrm{OU},t}
\left[
\exp\left(
\int_0^t
\frac{4(1-2Y_s^2)}{(1+2Y_s^2)^2}
\,\d s
\right)
\right]
=
1+
\frac{2(e^{2t}-1)-8xy\sinh t}
{(1+2x^2)(1+2y^2)}.
\]

We can obtain a corresponding identity for Brownian bridges
using the standard space-time transformation between Brownian
motion and the OU process (see \cite{BarczyKern2013,BarczyKern2013b} for
representations of OU bridges).
More precisely, fix $t>0$ and set
$T:=\frac12\log(1+2t)$.
Let $\{B_s\}_{0\leqslant s\leqslant t}$ be a Brownian bridge from
$x$ to $y$ over $[0,t]$ and $\mathbb E_{x\to y}^{\mathrm{BB},t}$ denote the expectation under this bridge law. Then
\[
Y_u:=e^{-u}B_{(e^{2u}-1)/2},
\qquad 0\leqslant u\leqslant T,
\]
is an OU bridge from $x$ to $y/\sqrt{1+2t}$
over $[0,T]$.
Applying the preceding identity and changing variables,
we obtain
\[
\mathbb E_{x\to y}^{\mathrm{BB},t}
\left[
\exp\left(
\int_0^t
\frac{4(1+2s-2B_s^2)}
{(1+2s+2B_s^2)^2}
\,\d s
\right)
\right]
=
1+
\frac{4t(1+2t-2xy)}
{(1+2x^2)(1+2t+2y^2)}.
\]

To introduce our next applications, for a strictly increasing sequence
$\k\in\mathbb N^\ell$, define
\[
\Phi_\k(\lambda):=
\frac{\Gamma((\lambda+1)/2)}{\Gamma(\lambda/2)}
\frac{\prod_{j:\,k_j\ {\rm odd}}(\lambda+k_j)}
{\prod_{j:\,k_j\ {\rm even}}(\lambda+k_j)}.
\]
The function $\Phi_\k$ extends meromorphically to
$\mathbb C$ and admits the following Weierstrass product:
\[
\Phi_\k(\lambda)
=
\frac{\sqrt{\pi}}{2}\lambda
\frac{\prod_{j:\,k_j\ {\rm odd}}k_j}
{\prod_{j:\,k_j\ {\rm even}}k_j}
\prod_{n=1}^\infty
\frac{
(1+\frac{\lambda}{2n}\mathbf 1_{\{2n\notin\k\}})
}{
(1+\frac{\lambda}{2n-1}\mathbf 1_{\{2n-1\notin\k\}})
}.
\]
This formula follows from the infinite product identity
\cite[Eq.~(5.8.5)]{NIST}.
We leave it as an exercise to check that the zeros and poles of $\Phi_\k$ interlace if and only if
$\k$ is a Krein--Adler sequence. Real meromorphic functions with interlacing zeros and poles are closely related to Pick functions (see  \cite{ChebotarevMeiman1949,Levin1996} or \cite[Theorem 2]{Kuznetsov_2010}). 

\begin{proposition}\label{prop:occupation-time}
Let $\k$ be a Krein--Adler sequence.
Then for $\lambda>0$,
\[
\Phi_\k(\lambda)
=
\frac12
\frac{(f_\lambda^+)'(0)}{f_\lambda^+(0)},
\]
where $f_\lambda^+$ is defined in
\eqref{def_f_lambda_pm}.
\end{proposition}

The logarithmic derivative of a fundamental
solution plays an important role in the study of local
times, occupation times and excursions of one-dimensional
diffusions. We now give several examples. For the next three applications, assume that $X_0^\k=0$, and define the local time
of $X^\k$ at zero by
\[
L_t^0:=
\lim_{\epsilon\downarrow0}
\frac{1}{2\epsilon}
\int_0^t
\mathbf 1_{\{|X_s^\k|<\epsilon\}}\,\d s.
\]
Let $\tau_a:=\inf\{t\geqslant 0:L_t^0>a\}$
be the inverse local time. By
\cite[Theorem~1, Eqs.~(14)--(17)]{PitmanYor2003},
$\{\tau_a\}_{a\geqslant 0}$ is a subordinator with Laplace
exponent $2\Phi_\k$:
\[
\mathbb E_0[e^{-\lambda\tau_a}]
=
e^{-2a\Phi_\k(\lambda)},
\qquad a,\lambda\geqslant 0.
\]

Next, define the occupation time of the positive half-line
by 
$A_t:=\int_0^t
\mathbf 1_{\{X_s^\k>0\}}\,\d s$.
The general formula of
\cite[Corollary~8, Eq.~(68)]{PitmanYor2003}
gives the joint Laplace transform
\[
\int_0^\infty e^{-\lambda t}
\mathbb E_0
\left[e^{-\alpha L_t^0-qA_t}\right]\,\d t
=
\frac{
\Phi_\k(\lambda)/\lambda+
\Phi_\k(\lambda+q)/(\lambda+q)
}{
\alpha+\Phi_\k(\lambda)+\Phi_\k(\lambda+q)
},
\]
where $\lambda>0$ and $\alpha,q\geqslant 0$.

There are also simple expressions involving the last
zero before an independent exponential time.
Let $T$ be independent of $X^\k$ and exponentially distributed with rate
$\lambda>0$. Define the last zero before $T$ as 
$g_T:=\sup\{s<T:X_s^\k=0\}$.
By \cite[Corollary~8, Eqs.~(66)--(67)]{PitmanYor2003},
the random variables $g_T$ and $T-g_T$ are independent,
and their Laplace transforms are
\[
\mathbb E_0[e^{-qg_T}]
=
\frac{\Phi_\k(\lambda)}{\Phi_\k(\lambda+q)}, \qquad 
\mathbb E_0[e^{-q(T-g_T)}]
=
\frac{\lambda\Phi_\k(\lambda+q)}
{(\lambda+q)\Phi_\k(\lambda)},
\qquad q\geqslant 0.
\]

Finally, consider the first hitting time of zero,
$h_0:=\inf\{t\geqslant 0:X_t^\k=0\}$.
Suppose that the initial distribution of $X^\k$ is
its stationary distribution $\pi_\k$ given in   \eqref{def_pi_k}.
Applying the general identity \cite[Eq.~(33)]{PitmanYor2003} to the stationary diffusion \(X^{\mathbf k}\), and using the fact that the inverse local time at zero has Laplace exponent \(2\Phi_{\mathbf k}\), we obtain
\[
\mathbb E_{\pi_\k}[e^{-\lambda h_0}]
=
\frac{\Phi_\k(\lambda)}
{\lambda\Phi_\k'(0)},
\qquad \lambda>0.
\]
In particular, if the Krein--Adler sequence $\k$ consists of consecutive integers
$\k=(1,2,\ldots,\ell)$, 
then
\[
\Phi_\k(\lambda)
=
\frac{\lambda}{2}
\frac{\Gamma((\lambda+\ell+1)/2)}
{\Gamma((\lambda+\ell+2)/2)}.
\]
Consequently, under the stationary initial distribution $X_0^{\k} \sim \pi_{\k}$, the random variable 
$\exp(-2h_0)$ has a $\mathrm{Beta}\left((\ell+1)/2,1/2 \right)$ distribution.

The paper is organized as follows.
In Section~\ref{Section_Q_n}, we study the connection polynomials
$Q_n^{\k}(x,y)$ and prove Theorem~\ref{thm:Mehler_kernel}.
We establish these results for arbitrary strictly increasing
sequences $\k$, as they are of independent interest in the
theory of Hermite Wronskians.
Section~\ref{Section_proofs} contains the proofs of all remaining
results.

%***************************************************************
%***************************************************************
%***************************************************************

\section{\texorpdfstring
  {Polynomials $Q_n^\k(x,y)$}
  {Connection polynomials}}
\label{Section_Q_n}

%***************************************************************
%***************************************************************
%***************************************************************

We begin by proving an identity for Wronskians of Hermite polynomials. 

\begin{lemma}\label{lemma:H_k,m}
Let $\k\in\n^\ell$ be strictly increasing, with $\ell\geqslant 1$.
For all $n\in\z_{\geqslant0}\setminus\k$ and $x,y\in\mathbb C$,
\begin{equation}\label{eq:Wr-partial-fractions}
\begin{aligned}
\frac{H_{\k,n}(x)H_{\k,n}(y)}
{\prod_{j=1}^{\ell}(n-k_j)}
={}&\mathbf 1_{\{\ell\ \mathrm{even}\}}2^\ell
H_{\k}(x)H_{\k}(y)H_n(x)H_n(y)\\
&+2^{2\ell-2}\sum_{i=1}^{\ell}
\frac{\prod_{j\ne i}(k_i-k_j)}{n-k_i}
H_{\k\setminus k_i}(x)H_{\k\setminus k_i}(y)
H_{k_i,n}(x)H_{k_i,n}(y).
\end{aligned}
\end{equation}
Moreover, for each $a=1,\ldots,\ell$, the right-hand side
of \eqref{eq:Wr-partial-fractions} vanishes at $n=k_a$, provided that the term $i=a$ is omitted in the finite sum.
\end{lemma}
\begin{proof}
Fix $n\geqslant 0$. We first rewrite the rows of $H_{\k,n}(x)$
using the Hermite differential operator
${\mathcal D}:=2\mathcal L^{\mathrm{OU}}=\partial_x^2-2x\partial_x$.
Since ${\mathcal D} H_m=-2mH_m$, we have
\[
{\mathcal D}^rH_m=(-2m)^rH_m,\qquad
\partial_x {\mathcal D}^rH_m=(-2m)^rH_m',
\qquad r\geqslant 0.
\]
The operators ${\mathcal D}^r$ and $\partial_x {\mathcal D}^r$ have leading
terms $\partial_x^{2r}$ and $\partial_x^{2r+1}$,
respectively. All remaining derivatives have lower
orders. Thus, replacing the operator $\partial_x^{2r}$ in row $2r$ by ${\mathcal D}^{r}$ (and the operator $\partial_x^{2r+1}$ in row $2r+1$ by $\partial_x {\mathcal D}^{r}$)  amounts to elementary row operations
and preserves the determinant. Numbering the rows from $0$ to $\ell$, we obtain
the following row pattern:
\begin{equation}\label{eq:general-row-pattern}
\begin{array}{c|cccc}
{\textnormal{ row $2r$}}
&(-2k_1)^rH_{k_1}&\cdots&(-2k_\ell)^rH_{k_\ell}
&(-2n)^rH_n\\[2mm]
{\textnormal{ row $2r+1$}}
&(-2k_1)^rH_{k_1}'&\cdots&(-2k_\ell)^rH_{k_\ell}'
&(-2n)^rH_n'
\end{array}
\end{equation}
where only rows with indices at most $\ell$ are included.
Here and below, we suppress the argument $x$ in
determinants.

Introduce a parameter $s$ by replacing $(-2n)^r$
in the last column by $(-2s)^r$, and denote the
resulting determinant by $P_{\k,s}(x)$.
Throughout this construction, $n$ is fixed:
the functions $H_n(x)$ and $H_n'(x)$ remain unchanged.
In particular, $P_{\k,n}(x)=H_{\k,n}(x)$. Only the last column depends on $s$, and its entries
have degree at most $\lfloor\ell/2\rfloor$.
Expansion along this column therefore shows that
$\deg_sP_{\k,s}(x)\leqslant \lfloor\ell/2\rfloor$.
If $\ell$ is even, only the last entry contributes
to the coefficient of $s^{\ell/2}$. Its cofactor
is $H_{\k}(x)$, so this coefficient is
$(-2)^{\ell/2}H_{\k}(x)H_n(x)$.

We next evaluate $P_{\k,k_i}(x)$. Perform the row
operations
\begin{equation}\label{row_operations}
{\textnormal{row $r$}}\longleftarrow ({\textnormal{row $r$}})+2k_i \times ({\textnormal{row $r-2$}}),
\qquad r=\ell,\ell-1,\ldots,2.
\end{equation}
The descending order ensures that each operation
uses the original row $r-2$. For an even row
$r=2q$, the entry in column $j$ becomes
\[
(-2k_j)^qH_{k_j}
+2k_i(-2k_j)^{q-1}H_{k_j}
=2(k_i-k_j)(-2k_j)^{q-1}H_{k_j}.
\]
The same calculation applies to odd rows, with
$H_{k_j}'$ in place of $H_{k_j}$. Since $s=k_i$,
the $i$th and last columns vanish below their
first two entries.

Expanding along these two columns gives the
factor $H_{k_i,n}(x)$ and the sign $(-1)^{\ell+i}$.
In the complementary determinant, we extract
$2(k_i-k_j)$ from each column $j\ne i$.
The remaining determinant has the row pattern
\eqref{eq:general-row-pattern} and equals
$H_{\k\setminus k_i}(x)$. Hence
\begin{equation}\label{P_k_k_i_identity}
P_{\k,k_i}(x)
=(-1)^{\ell+i}2^{\ell-1}
\Big[\prod_{j\ne i}(k_i-k_j)\Big]
H_{\k\setminus k_i}(x)H_{k_i,n}(x).
\end{equation}

Let us illustrate this calculation with $\ell=4$
and $i=2$. After the row operations
\eqref{row_operations}, the determinant becomes
\[
P_{\k,k_2}(x)=
\begin{vmatrix}
H_{k_1}&H_{k_2}&H_{k_3}&H_{k_4}&H_n\\
H_{k_1}'&H_{k_2}'&H_{k_3}'&H_{k_4}'&H_n'\\
2(k_2-k_1)H_{k_1}&0&2(k_2-k_3)H_{k_3}
&2(k_2-k_4)H_{k_4}&0\\
2(k_2-k_1)H_{k_1}'&0&2(k_2-k_3)H_{k_3}'
&2(k_2-k_4)H_{k_4}'&0\\
-4k_1(k_2-k_1)H_{k_1}&0&-4k_3(k_2-k_3)H_{k_3}
&-4k_4(k_2-k_4)H_{k_4}&0
\end{vmatrix}.
\]
We expand along columns $2$ and $5$. Their first
two entries give
\[
\begin{vmatrix}
H_{k_2}&H_n\\
H_{k_2}'&H_n'
\end{vmatrix}
=H_{k_2,n}(x).
\]
The cofactor sign is $(-1)^{4+2}=1$, so
\[
P_{\k,k_2}(x)=H_{k_2,n}(x)
\begin{vmatrix}
2(k_2-k_1)H_{k_1}&2(k_2-k_3)H_{k_3}
&2(k_2-k_4)H_{k_4}\\
2(k_2-k_1)H_{k_1}'&2(k_2-k_3)H_{k_3}'
&2(k_2-k_4)H_{k_4}'\\
-4k_1(k_2-k_1)H_{k_1}&-4k_3(k_2-k_3)H_{k_3}
&-4k_4(k_2-k_4)H_{k_4}
\end{vmatrix}.
\]
Extracting $2(k_2-k_1)$, $2(k_2-k_3)$ and
$2(k_2-k_4)$ from the three columns gives
\[
P_{\k,k_2}(x)
=8(k_2-k_1)(k_2-k_3)(k_2-k_4)H_{k_2,n}(x)
\begin{vmatrix}
H_{k_1}&H_{k_3}&H_{k_4}\\
H_{k_1}'&H_{k_3}'&H_{k_4}'\\
-2k_1H_{k_1}&-2k_3H_{k_3}&-2k_4H_{k_4}
\end{vmatrix}.
\]
Adding $2x$ times the second row to the third row
replaces each entry $-2k_jH_{k_j}$ by $H_{k_j}''$.
Thus the last determinant is
$H_{(k_1,k_3,k_4)}(x)=H_{\k\setminus k_2}(x)$,
and we obtain
\[
P_{\k,k_2}(x)
=8(k_2-k_1)(k_2-k_3)(k_2-k_4)
H_{\k\setminus k_2}(x)H_{k_2,n}(x),
\]
in agreement with \eqref{P_k_k_i_identity}.

We now use the degree bound and the values
\eqref{P_k_k_i_identity} to obtain a
partial-fraction identity, which would imply \eqref{eq:Wr-partial-fractions}. Set
$\pi(s):=\prod_{j=1}^{\ell}(s-k_j)$ and
$R(s):=P_{\k,s}(x)P_{\k,s}(y)/\pi(s)$.
The numerator of $R$ has degree at most $\ell$, and
the denominator has distinct zeros $k_1,\ldots,k_\ell$.
Consequently, the partial fraction decomposition of $R$ has the form
\begin{equation}\label{R_partial_fraction}
R(s)=A+\sum_{i=1}^{\ell}
\frac{P_{\k,k_i}(x)P_{\k,k_i}(y)}
{(s-k_i)\pi'(k_i)}.
\end{equation}
Here $A$ is the coefficient of $s^\ell$ in
$P_{\k,s}(x)P_{\k,s}(y)$. When $\ell$ is odd, the degree of the polynomial $s \mapsto P_{\k,s}(x)P_{\k,s}(y)$ is at most $2\lfloor \ell/2 \rfloor =\ell -1$, thus $A=0$ in this case. When $\ell$ is even, we find that $A=2^\ell H_{\k}(x)H_{\k}(y)H_n(x)H_n(y)$, since we established above that the coefficient of $s^{\ell/2}$ of the polynomial $s\mapsto P_{\k,s}(x)$ is $(-2)^{\ell/2}H_{\k}(x)H_n(x)$. Substituting
\eqref{P_k_k_i_identity}, using
$\pi'(k_i)=\prod_{j\ne i}(k_i-k_j)$, and setting
$s=n$ proves \eqref{eq:Wr-partial-fractions}
for $n\notin\k$.

Finally, suppose that $n=k_a$. At $s=k_a$, the
$a$th and last columns coincide. Thus both
$P_{\k,s}(x)$ and $P_{\k,s}(y)$ are divisible
by $s-k_a$. Their product has a zero of
multiplicity at least two, whereas $\pi(s)$
has a simple zero. It follows that
$R(s)\to0$ as $s\to k_a$. Moreover, the coefficient of the term $i=a$
in \eqref{R_partial_fraction} is zero.
Taking $s\to k_a$ in the remaining terms
proves the last assertion.
\end{proof}

As in \cite{Kuznetsov2026}, for $|z|<1$ and $x,y\in {\mathbb C}$ we denote
\begin{align}
E(z,x,y)&:=\sum_{n\geqslant0}
\frac{z^n}{2^n n!}H_n(x)H_n(y)=\frac{1}{\sqrt{1-z^2}} \exp\bigg( \frac{2xzy-z^2(x^2+y^2)}{(1-z^2)} \bigg),
\label{eq:E}\\
E_{\k}(z,x,y)&:=\sum_{n\in\z_{\geqslant0}\setminus\k}
\frac{z^n H_{\k,n}(x)H_{\k,n}(y)}
{2^{n+\ell}n!\prod_{j=1}^{\ell}(n-k_j)}.
\label{eq:Ek}
\end{align}
The explicit expression for $E(z,x,y)$ is the classical Mehler formula (see \cite[formula (18.18.28)]{NIST}). 

\vspace{0.25cm}
\noindent
\textbf{Proof of Theorem~\ref{thm:Mehler_kernel}:}
By Lemma~\ref{lemma:H_k,m}, the identity
\eqref{eq:Wr-partial-fractions} extends to all
$n\geqslant 0$ if its left-hand side is set to zero
for $n\in\k$ and the term $i=a$ is omitted
when $n=k_a$. Multiply this identity by $z^n/(2^{n+\ell}n!)$
and sum over $n\geqslant 0$. The series converge
absolutely for $|z|<1$ and $x,y\in\mathbb C$
\cite[Lemma~2.1(ii)]{Kuznetsov2026}, giving
\begin{equation}\label{eq:single-index-reduction}
\begin{aligned}
E_{\k}(z,x,y)
={}&\mathbf 1_{\{\ell\ \mathrm{even}\}}
H_{\k}(x)H_{\k}(y)E(z,x,y)\\
&+2^{\ell-1}\sum_{i=1}^{\ell}
\Big[\prod_{j\ne i}(k_i-k_j)\Big]
H_{\k\setminus k_i}(x)H_{\k\setminus k_i}(y)
E_{(k_i)}(z,x,y).
\end{aligned}
\end{equation}

The single-index case of the generalized Mehler
formula \cite[Theorem~1.1 and Proposition~1.5]{Kuznetsov2026}
gives
\[
E_{(k_i)}(z,x,y)
=E(z,x,y)\sum_{n=0}^{k_i+1}z^nQ_n^{(k_i)}(x,y),
\]
where, for $n\geqslant 1$,
\[
Q_n^{(k_i)}(x,y)
=2^{k_i}k_i!
\big[T_{k_i+1-n}(x,y)-T_{k_i-1-n}(x,y)\big].
\]
We set $Q^{(k_i)}_n=0$ for $n>k_i+1$, consistently
with this formula.

Substituting these identities into
\eqref{eq:single-index-reduction} shows that
$E_{\k}$ equals $E$ multiplied by a polynomial
in $z$ of degree at most $k_\ell+1$.
For $n\geqslant 1$, the coefficient of $z^n$ in this
polynomial is
\[
2^{\ell-1}\sum_{i=1}^{\ell}
\Big[\prod_{j\ne i}(k_i-k_j)\Big]
H_{\k\setminus k_i}(x)H_{\k\setminus k_i}(y)
Q_n^{(k_i)}(x,y).
\]
Comparing coefficients with \eqref{eqn:Mehler} and substituting
the preceding expression for $Q_n^{(k_i)}$ yields
\eqref{eq:explicit}.
\qed

\vspace{0.25cm}

Formula \eqref{eq:explicit} gives short proofs of
Conjectures~1.7 and~1.8 and the first assertion of
Conjecture~1.9 in \cite{Kuznetsov2026}.

Setting $n=k_\ell+1$, we observe that all terms
with $i<\ell$ in \eqref{eq:explicit} vanish.
Since $T_0(x,y)=1$ and $T_{-2}(x,y)=0$, we obtain
\[
Q^\k_{k_\ell+1}(x,y)
=2^{k_\ell+\ell-1}k_\ell!
\Big[\prod_{j=1}^{\ell-1}(k_\ell-k_j)\Big]
H_{\k\setminus k_\ell}(x)H_{\k\setminus k_\ell}(y),
\]
which proves Conjecture~1.7.

Next, take $\k=(k,k+1)$, with $k\geqslant 1$, and set
$n=k+1$. Using $T_0(x,y)=1$, $T_1(x,y)=2xy$,
and $T_m(x,y)=0$ for $m<0$, formula
\eqref{eq:explicit} gives
\[
Q^\k_{k+1}(x,y)
=2^{k+1}k!\big[4(k+1)xyH_k(x)H_k(y)
-H_{k+1}(x)H_{k+1}(y)\big].
\]
Applying the recurrence
$2xH_k(x)=H_{k+1}(x)+2kH_{k-1}(x)$
in both variables gives the expression stated
in Conjecture~1.8.

To address Conjecture~1.9, set $h_m(x):=H_m(x/2)$.
The recurrence $h_{m+1}=xh_m-2mh_{m-1}$, together
with $h_0=1$ and $h_1=x$, shows by induction that
$h_m\in\mathbb Z[x]$ for all $m\geqslant 0$.
For a sequence $\j=(j_1,\ldots,j_r)$, the chain
rule gives
\[
H_\j(x/2)
=2^{r(r-1)/2}\Wr[h_{j_1},\ldots,h_{j_r}](x)
\in\mathbb Z[x].
\]
Moreover, when $n\geqslant 1$, every nonzero $T_m$ term
in the $i$th summand of \eqref{eq:explicit}
has $0\leqslant m\leqslant k_i$. Hence
\[
2^{k_i}k_i!T_m(x/2,y/2)
=2^{k_i-m}\frac{k_i!}{m!}h_m(x)h_m(y)
\in\mathbb Z[x,y].
\]
Formula \eqref{eq:explicit} therefore gives
$Q_n^\k(x/2,y/2)\in\mathbb Z[x,y]$
for $1\leqslant n\leqslant k_\ell+1$.
The case $n=0$ follows directly from
\eqref{eqn_Qk0}. 
Writing 
\[\Delta(\k):=\prod_{1\leqslant i<j \leqslant \ell}(k_j-k_i),\]
the stronger assertion
$\Delta(\k)^{-2}Q_n^\k(x/2,y/2)\in\mathbb Z[x,y]$
in Conjecture~1.9 is not established by this argument.

%***************************************************************
%***************************************************************
%***************************************************************
\section{Proofs} \label{Section_proofs}
%***************************************************************
%***************************************************************
%***************************************************************

Throughout this section, $\k$ denotes a Krein--Adler sequence.
We write
\begin{equation}\label{def_|k|_h_k}
|\k|:=\sum_{i=1}^{\ell}k_i,\quad
h_{\k}(x):=e^{-x^2/2}\frac{H_{\k-1}(x)}{H_{\k}(x)},
\end{equation}
and define
\[
\mathcal H^{\k}
:=-\partial_x^2+x^2-2(\log H_{\k})''.
\]

\begin{lemma}\label{Lemma_H_k_H_k-1}
Let $\k$ be a Krein--Adler sequence.
The polynomials $H_{\k}$ and $H_{\k-1}$ are even
and strictly positive on $\r$. Moreover, $H_{\k}$
has degree $|\k|-\ell(\ell-1)/2$ and leading
coefficient $2^{|\k|}\Delta(\k)$.
\end{lemma}

\begin{proof}
Let $\j=(j_1,\ldots,j_\ell)$ be a strictly increasing
sequence of nonnegative integers.
The standard degree formula for Wronskians
\cite[Lemma~2.1]{BonneuxStevens2018} shows that
$H_\j$ has degree $|\j|-\ell(\ell-1)/2$
and leading coefficient
$2^{|\j|}\prod_{1\leqslant a<b\leqslant \ell}(j_b-j_a)>0$.
This proves the degree and leading-coefficient
assertions for $H_\k$ and also shows that
$H_{\k-1}$ has a positive leading coefficient.

Recall that Hermite functions are defined by $\psi_n(x):=\exp(-x^2/2) H_n(x)$. By the Krein--Adler theorem (\cite{Adler1994} or 
\cite[Theorem~4.2]{Gomez-Ullate_2014}),
the Wronskians of Hermite functions indexed
by $\k$ and $\k-1$ have no real zeros.
The theorem allows an initial block beginning
at $0$, as may occur for $\k-1$.
Since
$\Wr[\psi_{j_1},\ldots,\psi_{j_\ell}](x)
=e^{-\ell x^2/2}H_\j(x)$,
the same holds for $H_\k$ and $H_{\k-1}$.
Their positive leading coefficients therefore
imply that both polynomials are strictly
positive on $\r$.

Finally, \cite[Lemma~3.6]{BonneuxStevens2018} gives
\[
H_\j(-x)=(-1)^{|\j|-\ell(\ell-1)/2}H_\j(x).
\]
Thus each polynomial is either even or odd.
Since neither $H_\k$ nor $H_{\k-1}$ vanishes
at $0$, both are even.
\end{proof}

\vspace{0.25cm}
\noindent
\textbf{Proof of Theorem~\ref{thm:spectrum}.}
 Multiplication by $H_{\k-1}$ is a unitary map from
$L_2(\r,\nu_\k(x)\,\d x)$ onto \\
$L_2(\r,e^{-x^2}H_\k(x)^{-2}\,\d x)$,
sending $\varphi_{\k,n}$ to $H_{\k,n}$.
Thus orthogonality and the norm formula follow from
\cite[Corollary~6.3 and Theorem~6.4]{Duran_2014},
while completeness follows from the corrected proof
in \cite[Theorem~2.1]{Gomez-Ullate_2020}.

To obtain the eigenfunction equation, first apply
\cite[Proposition~5.2]{Gomez-Ullate_2014} with $n=0$,
using the identity
\begin{equation}\label{H_k0_identity}
H_{\k,0}=2^\ell\Big[\prod_{j=1}^{\ell}k_j\Big]H_{\k-1},
\end{equation}
which follows from \cite[Lemma~3.1]{Ullate_2018},
since $\ell$ is even.
Combining the resulting differential equation for
$H_{\k-1}$ with a direct calculation gives
\[
H_{\k-1}\mathcal L^\k\frac{1}{H_{\k-1}}
=\frac12\partial_x^2
-\Big(x+\frac{H_\k'}{H_\k}\Big)\partial_x
+\frac12\frac{H_\k''}{H_\k}
+x\frac{H_\k'}{H_\k}-\ell.
\]
Applying the same proposition to $H_{\k,n}$ now yields
$\mathcal L^\k\varphi_{\k,n}=-n\varphi_{\k,n}$,
as required. \qed

\vspace{0.25cm}
\noindent
\textbf{Proof of Theorem~\ref{thm:classification}.}
For $R\equiv0$, solvability follows from the classical
Hermite polynomials; for $R$ of the form \eqref{def_R_k},
it follows from Theorem~\ref{thm:spectrum}.
We therefore prove the converse.

Suppose that $\mathcal L_R$ (defined in \eqref{def_L_R}) is solvable by polynomials,
and set $h:=\sqrt{\mu_R}$.
Since $h'/h=-(x+R)$, a direct calculation gives
\[
\mathcal H_R:=-2h\mathcal L_Rh^{-1}
=-\partial_x^2+(x+R)^2-1-R'
=-\partial_x^2+\frac{h''}{h}.
\]
Set $C:=\lim_{x\to\infty}(2xR(x)-1)$.
This limit exists because $R$ is rational and
$R(x)=O(|x|^{-1})$.
The operator $\widetilde{\mathcal H}:=\mathcal H_R-C$
has potential $x^2+\widetilde R(x)$, where
$\widetilde R=2xR+R^2-R'-1-C$ is rational, has no
real poles, and tends to zero as $|x|\to\infty$.
Thus $\widetilde{\mathcal H}$ is a regular rational
extension of the harmonic oscillator in the sense of
\cite[Definition~1.2]{Gomez-Ullate_2014}.

Multiplication by $h$ is a unitary map from
$L_2(\r,\mu_R(x)\,\d x)$ onto $L_2(\r,\d x)$
and takes each eigenfunction $wp_n$ to $(hw)p_n$.
Consequently, $\widetilde{\mathcal H}$ is solvable
by polynomials.
By \cite[Theorem~1.1]{Gomez-Ullate_2014}, either
$\widetilde{\mathcal H}=\mathcal H$ or
$\widetilde{\mathcal H}=\mathcal H^\k$
for some Krein--Adler sequence $\k$.

Since $\mathcal H_Rh=0$, the function $h$ is a strictly
positive $L_2$ eigenfunction of $\widetilde{\mathcal H}$
with eigenvalue $-C$.
It must therefore be a ground-state eigenfunction,
which is unique up to a constant factor
\cite[Theorem~10.13]{Teschl2014}.
If $\widetilde{\mathcal H}=\mathcal H$, this gives
$h(x)=c e^{-x^2/2}$ and hence $R\equiv0$. Now suppose that $\widetilde{\mathcal H}=\mathcal H^\k$,
where $\k$ has length $\ell$.
By \eqref{H_k0_identity}, $h_\k$ (defined in \eqref{def_|k|_h_k}) is proportional to
$e^{-x^2/2}H_{\k,0}/H_\k$.
The eigenfunction equation in the proof of
\cite[Proposition~5.5]{Gomez-Ullate_2014} therefore gives
\begin{equation}\label{eq:ground-state-Hk}
\mathcal H^\k h_\k=(1-2\ell)h_\k.
\end{equation}
Lemma~\ref{Lemma_H_k_H_k-1} and the degree formulas
show that $h_\k$ is strictly positive and belongs to
$L_2(\r,\d x)$.
Ground-state uniqueness now implies that
$C=2\ell-1$ and $h=c\,h_\k$ for some $c>0$.
Taking logarithmic derivatives, we obtain
\[
-x-R=\frac{h'}{h}
=\frac{h_\k'}{h_\k}
=-x+\frac{H_{\k-1}'}{H_{\k-1}}-\frac{H_\k'}{H_\k},
\]
which proves \eqref{def_R_k}.
\qed

\vspace{0.25cm}
We record the resulting conjugation identity
\begin{equation}\label{eq:Lk-conjugation}
\mathcal L^\k
=-\frac12 h_\k^{-1}
\big(\mathcal H^\k-(1-2\ell)\big)h_\k,
\end{equation}
which will be used in the proof of
Theorem~\ref{thm:rational-ou-fundamental}.

\vspace{0.25cm}
\noindent
\textbf{Proof of Theorem~\ref{thm:transition-density}.}
The transition probability density $p_t^\k(x,y)$ is the integral kernel of the operator $\exp(t \mathcal L^\k)$ \cite{Avramidi2015,Vassilevich2003}.  
By Theorem~\ref{thm:spectrum}, it has spectral expansion
\[
p_t^\k(x,y)
=\nu_\k(y)
\sum_{n\in\z_{\geqslant0}\setminus\k}
e^{-nt}
\frac{\varphi_{\k,n}(x)\varphi_{\k,n}(y)}
{\lVert\varphi_{\k,n}\rVert^2},
\]
where the norm is taken in $L_2(\r,\nu_\k(x)\,\d x)$.
The series converges absolutely and uniformly on
compact subsets of $(0,\infty)\times\r^2$,
by \cite[Lemma~2.1(ii)]{Kuznetsov2026}.

Substituting $\varphi_{\k,n}=H_{\k,n}/H_{\k-1}$
and the norm formula \eqref{eq:norm}, and using
the definition \eqref{eq:Ek}, we obtain
\[
p_t^\k(x,y)
=\frac{e^{-y^2}}{\sqrt\pi}
\frac{H_{\k-1}(y)}
{H_\k(y)^2H_{\k-1}(x)}
E_\k(e^{-t},x,y).
\]
On the other hand, \eqref{p_t^OU_formula} and
\eqref{eq:E} give
$p_t^{\mathrm{OU}}(x,y)
=e^{-y^2}E(e^{-t},x,y)/\sqrt\pi$.
Applying \eqref{eqn:Mehler} with $z=e^{-t}$
now yields \eqref{eqn:p_t^X}.
\qed

\vspace{0.25cm}
\noindent
\textbf{Proof of Theorem~\ref{thm:rational-ou-fundamental}.}
Fix $\lambda>0$ and write
$f:=F_{\k,\lambda}/H_{\k-1}$.
Set $u_\lambda:=e^{-x^2/2}F_\lambda$.
Recall that $\mathcal H=-\partial_x^2+x^2$ and
$\psi_n=e^{-x^2/2}H_n$. We have
\begin{equation}\label{eqn:H_psi_n}
\mathcal H u_\lambda=(1-2\lambda)u_\lambda,
\qquad
\mathcal H\psi_n=(2n+1)\psi_n.
\end{equation}
Since
$\Wr[\psi_{k_1},\ldots,\psi_{k_\ell}]
=e^{-\ell x^2/2}H_\k$,
the Darboux--Crum transformation \cite{Crum_1955}
maps $\mathcal H$ to $\mathcal H^\k+2\ell$.
Consequently, the function
\begin{equation}\label{eq:Uklambda}
U_{\k,\lambda}
:=\frac{\Wr[\psi_{k_1},\ldots,\psi_{k_\ell},u_\lambda]}
{\Wr[\psi_{k_1},\ldots,\psi_{k_\ell}]}
=e^{-x^2/2}\frac{F_{\k,\lambda}}{H_\k}
\end{equation}
satisfies
$\mathcal H^\k U_{\k,\lambda}
=(1-2\ell-2\lambda)U_{\k,\lambda}$.
The Crum identities initially hold away from zeros
of the intermediate Wronskians and extend to all real
$x$ by continuity, since the final denominator $H_\k$
has no real zeros.
Since $f=U_{\k,\lambda}/h_\k$,
\eqref{eq:Lk-conjugation} gives
$\mathcal L^\k f=\lambda f$.

We next examine the behavior of $f(x)$ as $x\to -\infty$.
The standard parabolic-cylinder formulas
\cite[Eqs.~(12.2.5), (12.8.9), and (12.9.1)]{NIST}
give, for each fixed integer $r\geqslant 0$,
\begin{equation}\label{eq:F-derivative}
F_\lambda^{(r)}(-t)
=2^{-\lambda/2}(\lambda)_r
t^{-\lambda-r}(1+O(t^{-2})),
\qquad t\to\infty,
\end{equation}
where $(\lambda)_r:=\Gamma(\lambda+r)/\Gamma(\lambda)$. We will show that, as $t\to\infty$,
\begin{equation}\label{eq:Fklambda-asymptotic}
F_{\k,\lambda}(-t)
=
2^{|\k|-\lambda/2}\Delta(\k)
\prod_{j=1}^{\ell}(\lambda+k_j)
t^{|\k|-\ell(\ell+1)/2-\lambda}(1+O(t^{-2})).
\end{equation}

Consider the Wronskian determinant representing
$F_{\k,\lambda}(-t)$. The leading term and parity of $H_m$ give
\[
H_m^{(r)}(-t)
=
2^m(-1)^{m-r}t^{m-r}
\big(m^{\underline r}+O(t^{-2})\big),
\]
where
$a^{\underline r}:=a(a-1)\cdots(a-r+1)$
is the falling factorial, with $a^{\underline0}:=1$.
Using the preceding formula and \eqref{eq:F-derivative},
we factor $2^{k_j}(-1)^{k_j}t^{k_j}$ from column $j$,
for $1\leqslant j\leqslant \ell$, and $2^{-\lambda/2}t^{-\lambda}$
from the last column.
Numbering the rows by $r=0,\ldots,\ell$, we also factor
$(-1)^rt^{-r}$ from row $r$.
Since $(-1)^r(\lambda)_r=(-\lambda)^{\underline r}$,
this gives
\[
F_{\k,\lambda}(-t)
=
(-1)^{|\k|+d}2^{|\k|-\lambda/2}
t^{|\k|-d-\lambda}
\big(\det M+O(t^{-2})\big),
\]
where $d:=\ell(\ell+1)/2$ and
\[
M:=\left(a_j^{\underline r}\right)_{
0\leqslant r\leqslant \ell,\;1\leqslant j\leqslant \ell+1},
\qquad
(a_1,\ldots,a_{\ell+1})
:=(k_1,\ldots,k_\ell,-\lambda).
\]

Since $a^{\underline r}$ is a monic polynomial
of degree $r$, elementary row operations reduce $M$
to a Vandermonde matrix. Hence
\[
\det M
=\prod_{1\leqslant i<j\leqslant \ell+1}(a_j-a_i)
=(-1)^\ell\Delta(\k)
\prod_{j=1}^{\ell}(\lambda+k_j).
\]
Writing $\ell=2m$, the Krein--Adler condition gives
$|\k|\equiv m\pmod2$, while
$d=m(2m+1)\equiv m\pmod2$.
Thus $(-1)^{|\k|+d}=(-1)^\ell=1$, proving
\eqref{eq:Fklambda-asymptotic}.

Combining \eqref{eq:Fklambda-asymptotic} with the
degree and leading coefficient of $H_{\k-1}$ gives
\[
f(-t)
=
2^{\ell-\lambda/2}
\prod_{j=1}^{\ell}(\lambda+k_j)
t^{-\lambda}(1+O(t^{-2})).
\]
Thus $f$ is positive for sufficiently negative $x$
and tends to zero as $x\to-\infty$.

Suppose that $x_0$ is the first zero of $f$.
Then $f>0$ on $(-\infty,x_0)$, and the limiting value
$f(-\infty)=0$ 
 implies that $f$ attains a positive
maximum at some point $a<x_0$.
At this point, $f'(a)=0$ and $f''(a)\leqslant 0$.
However, $\mathcal L^\k f=\lambda f$ gives
$f''(a)=2\lambda f(a)>0$, a contradiction.
Hence $f>0$ on $\r$.

If $f'$ were negative somewhere, the limit
$f(-\infty)=0$ would likewise force a positive
local maximum, giving the same contradiction.
Thus $f'\geqslant 0$.
If $f'(a)=0$ at some point, the equation $\mathcal L^\k f=\lambda f$ gives
$f''(a)=2\lambda f(a)>0$, so $f'$ would be negative
immediately to the left of $a$.
Therefore $f'>0$ on $\r$, and we may take
$f_\lambda^+=f$ as the increasing fundamental solution.

The drift term in $\mathcal L^\k$ is odd (due to Lemma \ref{Lemma_H_k_H_k-1}), so $f(-x)$
also solves $\mathcal L^\k f=\lambda f$.
This function is positive and strictly decreasing; hence we can take
$f_\lambda^-(x)=f(-x)$ as the decreasing fundamental
solution.

It remains to compute the Wronskian of $f^{\pm}_{\lambda}(x)$.
We first record how a single Darboux transformation
changes a Wronskian.
Suppose that for some potential function $V$ the functions $h$, $u$ and $v$ satisfy
\[
-h''+Vh=\alpha h,
\qquad
-u''+Vu=\beta u,
\qquad
-v''+Vv=\beta v.
\]
For $\mathcal D^h u:=u'-(h'/h)u$, a direct calculation
gives
\begin{equation}\label{eq:Darboux-Wronskian}
\Wr[\mathcal D^h u,\mathcal D^h v]
=(\beta-\alpha)\Wr[u,v]
\end{equation}
on every interval where $h$ is nonzero.

We now specify the successive Darboux transformations.
Set $\Psi_0:=1$ and
\[
\Psi_j:=\Wr[\psi_{k_1},\ldots,\psi_{k_j}],
\qquad 1\leqslant j\leqslant \ell.
\]
Work initially on an interval where all these Wronskians
are nonzero, and define
\[
\mathcal K^{(j)}
:=-\partial_x^2+x^2-2(\log|\Psi_j|)'',
\qquad 0\leqslant j\leqslant \ell.
\]
Thus $\mathcal K^{(0)}=\mathcal H$.
Moreover, $\Psi_\ell=e^{-\ell x^2/2}H_\k$, so
$\mathcal K^{(\ell)}=\mathcal H^\k+2\ell$.

By the Darboux--Crum formulas
\cite[Sections~2--3]{Crum_1955},
the $j$th Darboux transformation maps
$\mathcal K^{(j-1)}$ to $\mathcal K^{(j)}$.
Its seed function and first-order operator are
\[
\eta_j:=\frac{\Psi_j}{\Psi_{j-1}},
\qquad
\mathcal D^{(j)}:=\partial_x-\frac{\eta_j'}{\eta_j}.
\]
Here $\eta_j$ is the image of $\psi_{k_j}$ under
the preceding $j-1$ transformations.
For example, $\eta_1=\psi_{k_1}$, whereas
\[
\eta_2
=\mathcal D^{(1)}\psi_{k_2}
=\frac{\Wr[\psi_{k_1},\psi_{k_2}]}{\psi_{k_1}}.
\]
The Darboux--Crum identities \cite{Crum_1955} give
\[
\mathcal K^{(j-1)}\eta_j=(2k_j+1)\eta_j,
\qquad
\mathcal K^{(j)}\mathcal D^{(j)}
=\mathcal D^{(j)}\mathcal K^{(j-1)}.
\]

We apply these transformations to the two solutions
$u_0^+(x):=u_\lambda(x)$ and
$u_0^-(x):=u_\lambda(-x)$.
Define recursively
$u_j^\pm:=\mathcal D^{(j)}u_{j-1}^\pm$.
Crum's formula gives the explicit expressions
\[
u_j^\pm
=\frac{1}{\Psi_j}\Wr[\psi_{k_1},\ldots,\psi_{k_j},u_0^\pm],
\qquad 1\leqslant j\leqslant \ell.
\]
The intertwining identity above and
\eqref{eqn:H_psi_n} imply
$\mathcal K^{(j)}u_j^\pm=(1-2\lambda)u_j^\pm$.

Thus, at step $j$, the three functions
$\eta_j,u_{j-1}^-,u_{j-1}^+$ solve equations
with the same operator $\mathcal K^{(j-1)}$.
We may therefore apply \eqref{eq:Darboux-Wronskian}
with
\[
h=\eta_j,\qquad
u=u_{j-1}^-,\qquad
v=u_{j-1}^+,\qquad
\alpha=2k_j+1,\qquad
\beta=1-2\lambda.
\]
Since $\mathcal D^{(j)}u_{j-1}^\pm=u_j^\pm$,
this gives
\[
\Wr[u_j^-,u_j^+]
=-2(\lambda+k_j)\Wr[u_{j-1}^-,u_{j-1}^+].
\]

At the final step, \eqref{eq:Uklambda} gives
$u_\ell^+=U_{\k,\lambda}=h_\k f_\lambda^+$.
The parity of the Hermite functions gives
$u_j^-(x)=(-1)^j u_j^+(-x)$.
Since $\ell$ is even and $h_\k$ is even, we obtain
$u_\ell^-(x)=u_\ell^+(-x)=h_\k(x)f_\lambda^-(x)$.
Iterating the preceding Wronskian identity therefore
yields
\begin{equation}\label{eq:Wr-Darboux-product}
\Wr[h_\k f_\lambda^-,h_\k f_\lambda^+]
=
(-2)^\ell
\Big[\prod_{j=1}^{\ell}(\lambda+k_j)\Big]
\Wr[u_\lambda(-\cdot),u_\lambda].
\end{equation}

Each intermediate Wronskian $\Psi_j$ is a nonzero polynomial
times a Gaussian, so its real zeros form a finite set.
The final transformed functions are smooth on $\r$
because $\Psi_\ell$ has no real zeros.
Consequently, \eqref{eq:Wr-Darboux-product} extends
to all real $x$ by continuity.

Finally, $u_\lambda(\pm x)=e^{-x^2/2}g_\lambda^\pm(x)$,
so \eqref{g^pm_lambda_Wronskian} gives
\[
\Wr[u_\lambda(-\cdot),u_\lambda]
=e^{-x^2}\Wr[g_\lambda^-,g_\lambda^+]
=\frac{2\sqrt\pi}{\Gamma(\lambda)}.
\]
Combining this with \eqref{eq:Wr-Darboux-product},
and using the evenness of $\ell$, we obtain
\[
h_\k^2\Wr[f_\lambda^-,f_\lambda^+]
=
\frac{2^{\ell+1}\sqrt\pi}{\Gamma(\lambda)}
\prod_{j=1}^{\ell}(\lambda+k_j).
\]
Since $h_\k^2=\nu_\k$, this proves
\eqref{f^pm_lambda_Wronskian}.
\qed

\vspace{0.25cm}
\noindent
{\bf Proof of Theorem \ref{thm:OU-Darb}:}
Write $N:=k_\ell$ and set
\[
A_\k:=\{a_1,\dots,a_p\}=\{N-r:r\in\{0,\dots,N\}\setminus\k\},\qquad p=N+1-\ell,
\]
where the sequence $a_j$ is increasing.
For $n\geqslant 0$, define
$\widetilde H_n(x):=\i^{-n}H_n(\i x)$ and
$\eta_n(x):=e^{x^2/2}\widetilde H_n(x)
=\i^{-n}\psi_n(\i x)$.
By \eqref{eqn:H_psi_n}, these functions satisfy
$\mathcal H\eta_n=-(2n+1)\eta_n$.

For $1\leqslant j\leqslant p$, put
$W_j:=\Wr[\eta_{a_1},\ldots,\eta_{a_j}]$ and
$D_j:=\{0,\ldots,a_j\}\setminus\{a_j-a_i:1\leqslant i\leqslant j\}$.
We also set $W_0=1$ and $D_0=\varnothing$,
and list all finite sets of indices in increasing order.
Since $\k$ consists of adjacent pairs, $a_1$ is even
and every difference $a_j-a_{j-1}$ is odd for $2\leqslant j\leqslant p$.
Thus $D_1=\{1,2,\ldots,a_1\}$ is a Krein--Adler sequence.
For $j\geqslant 2$, writing $d=a_j-a_{j-1}$, we have
$D_j=\{1,\ldots,d-1\}\cup(d+D_{j-1})$.
Since $d$ is odd, induction shows that every $D_j$
is a Krein--Adler sequence. Moreover, $a_p=N$
and $D_p=\k$.

The pseudo-Wronskian identities
\cite[Theorem~3.1 and Proposition~3.2]{Ullate_2018} give
\[
W_j(x)=B_j e^{jx^2/2}H_{D_j}(x).
\]
Here $B_j>0$, since
$\Wr[\widetilde H_{a_1},\ldots,\widetilde H_{a_j}]$
and $H_{D_j}$ have positive leading coefficients.
Lemma~\ref{Lemma_H_k_H_k-1} therefore implies that
$W_j$ is strictly positive on $\r$. Consequently,
\[
\rho_j(x):=\frac{W_j(x)}{W_{j-1}(x)}
=C_j e^{x^2/2}\frac{H_{D_j}(x)}{H_{D_{j-1}}(x)},
\qquad C_j>0,
\]
is smooth and strictly positive.

Define
\begin{equation}\label{eq:Hj-Darboux}
\mathcal H^{(j)}:=\mathcal H-2(\log W_j(x))''
=-\partial_x^2+x^2-2j-2(\log H_{D_j}(x))'',
\qquad 0\leqslant j\leqslant p.
\end{equation}
Let $Y^{(0)}$ be the Brownian motion on $\r$ killed at
rate $c_0(x)=x^2/2$. Starting from $Y^{(0)}$, we construct $Y^{(1)},\ldots,Y^{(p)}$
by successive Darboux transformations of \cite[Section~4]{KuznetsovYuan} with seeds $\rho_j$.
We prove inductively that their generators have the form
\begin{equation}\label{eq:Lj-Darboux}
\mathcal L^{(j)}
=\frac12\partial_x^2-c_j(x)
=-\frac12\mathcal H^{(j)}-\gamma_j,
\end{equation}
where $c_j(x) \geqslant 0$ for all $x\in \r$ and $\gamma_j$ is a real constant (and $\gamma_0=0$).

Suppose that $Y^{(j-1)}$ has been constructed.
By Crum's formula \cite{Crum_1955},
$\mathcal H^{(j-1)}\rho_j=-(2a_j+1)\rho_j$, and hence
\begin{equation}\label{eq:rho-eigenvalue}
\mathcal L^{(j-1)}\rho_j=\lambda_j\rho_j,
\qquad \lambda_j:=a_j+\frac12-\gamma_{j-1}.
\end{equation}
Set
\[
b_j(x):=\frac{\rho_j'(x)}{\rho_j(x)}
=x+\frac{H_{D_j}'(x)}{H_{D_j}(x)}
-\frac{H_{D_{j-1}}'(x)}{H_{D_{j-1}}(x)}.
\]
This function is smooth on $\r$ and satisfies
$b_j(x)=x+O(|x|^{-1})$ and
$b_j'(x)=1+O(|x|^{-2})$ as $|x|\to\infty$.
In particular, $b_j$ is globally Lipschitz.

We next verify that $\rho_j$ is $\lambda_j$-invariant.
For Brownian motion $B$ started at $x$, define
\[
M_t^{(j)}
:=\frac{\rho_j(B_t)}{\rho_j(x)}
\exp\left(-\lambda_jt-\int_0^t c_{j-1}(B_s)\,\d s\right).
\]
It\^o's formula and \eqref{eq:rho-eigenvalue} give
$\d M_t^{(j)}=M_t^{(j)}b_j(B_t)\,\d B_t$,
so $M^{(j)}$ is the stochastic exponential of
$\int_0^\cdot b_j(B_s)\,\d B_s$.
Since $b_j$ is globally Lipschitz, the auxiliary diffusion
$\d Z_t=b_j(Z_t)\,\d t+\d W_t$ is nonexplosive.
Brownian motion is also nonexplosive, so
\cite[Corollary~2.2]{MijatovicUrusov2012}
implies that $M^{(j)}$ is a martingale. Thus ${\mathbb E}_x[M^{(j)}_t]=1$, and rewriting this identity in terms of killed Brownian motion $Y^{(j-1)}$ gives 
\[
\mathbb E_x[\rho_j(Y_t^{(j-1)})]
=e^{\lambda_jt}\rho_j(x),\qquad t\geqslant 0,
\]
where $\rho_j$ vanishes at the cemetery state. Thus, the function $\rho_j$ is $\lambda_j$-invariant for $Y^{(j-1)}$. 

By \eqref{eq:Hj-Darboux} and \eqref{eq:Lj-Darboux},
$c_{j-1}(x)=x^2/2+O(1)$ as $|x|\to\infty$.
Thus $c_{j-1}(x)-b_j(x)^2=-x^2/2+O(1)$, and
\[
\mathfrak m_j
:=\sup_{x\in\r}\{c_{j-1}(x)-b_j(x)^2\}<\infty.
\]
This verifies Assumption~4.1 of \cite{KuznetsovYuan}.
The Doob transform of $Y^{(j-1)}$ by the martingale $M^{(j)}$ is a conservative diffusion on $\r$ with drift $b_j(x)$,
so the boundaries $\pm \infty$ are inaccessible by the nonexplosion
argument above. Hence Assumption~4.3 of  \cite{KuznetsovYuan} is also satisfied.

We may therefore apply the Darboux transformation of
\cite[Section~4]{KuznetsovYuan}.
By equation~(4.11) there, the resulting process $Y^{(j)}$
is Brownian motion killed at rate
\begin{equation}\label{eq:cj-Darboux}
c_j(x)=\mathfrak m_j+b_j(x)^2-c_{j-1}(x)\geqslant 0.
\end{equation}
The eigenfunction equation gives
$b_j'(x)=2(c_{j-1}(x)+\lambda_j)-b_j(x)^2$, while
$\mathcal H^{(j)}=\mathcal H^{(j-1)}-2b_j'(x)$.
Combining these identities with \eqref{eq:cj-Darboux}
yields
\[
\mathcal L^{(j)}=-\frac12\mathcal H^{(j)}-\gamma_j,
\qquad
\gamma_j:=\gamma_{j-1}+\mathfrak m_j+2\lambda_j.
\]
This completes the induction.

Since $D_p=\k$, we have
$\mathcal H^{(p)}=\mathcal H^\k-2p$, and therefore
\[
\mathcal L^{(p)}=-\frac12\mathcal H^\k-C,
\qquad C:=\gamma_p-p.
\]
By \eqref{eq:ground-state-Hk}, the bounded positive function
$h_\k$ satisfies
$\mathcal L^{(p)}h_\k=\lambda_*h_\k$, where
$\lambda_*:=\ell-\frac12-C$.
Since $c_p\geqslant 0$ and $h_\k$ is bounded, the local martingale
\[
M_t:=e^{-\lambda_*t-\int_0^t c_p(B_s)\,\d s}h_\k(B_t)
\]
is bounded on every finite time interval and is therefore a martingale.
We conclude that $h_\k$ is $\lambda_*$-invariant for $Y^{(p)}$
and defines a conservative Doob transform. By \eqref{eq:Lk-conjugation}, the generator of this
Doob transform is
\[
h_\k^{-1}(\mathcal L^{(p)}-\lambda_*)h_\k
=-\frac12h_\k^{-1}\bigl(\mathcal H^\k-(1-2\ell)\bigr)h_\k
=\mathcal L^\k.
\]
Hence the resulting diffusion is $X^\k$, as required.
\qed

\vspace{0.25cm}
\noindent
\textbf{Proof of Proposition~\ref{cor:OU-bridge}.}  Let
\[
\rho(x):=\frac{H_{\k-1}(x)}{H_\k(x)}
=e^{x^2/2}h_\k(x),
\qquad
V(x):=(\log H_\k)''(x).
\]
Recall that
$\mathcal L^{\mathrm{OU}}=\tfrac12\partial_x^2-x\partial_x$.
Substituting $h_\k=e^{-x^2/2}\rho$ into
\eqref{eq:ground-state-Hk} gives
\begin{equation}\label{eq:rho-eigenfunction}
(\mathcal L^{\mathrm{OU}}+V)\rho=\ell\rho.
\end{equation}

Let $Y$ be the OU process started at $x$, and define
\[
M_t:=
\frac{\rho(Y_t)}{\rho(x)}
\exp\left(-\ell t+\int_0^t V(Y_s)\,\d s\right).
\]
By It\^o's formula and \eqref{eq:rho-eigenfunction},
\[
\d M_t=M_t\theta(Y_t)\,\d W_t,
\qquad
\theta(x):=\frac{\rho'(x)}{\rho(x)}=-R_\k(x),
\]
where $R_\k$ denotes the right-hand side of
\eqref{def_R_k}.
Since $M_0=1$, it follows that
\[
M_t=
\exp\left(
\int_0^t\theta(Y_s)\,\d W_s
-\frac12\int_0^t\theta(Y_s)^2\,\d s
\right).
\]
The function $\theta$ is bounded, since $H_\k$ and
$H_{\k-1}$ have no real zeros and
$R_\k(x)=O(|x|^{-1})$ as $|x|\to\infty$.
Thus Novikov's condition holds on every finite time
interval, and $\{M_t\}_{t\geqslant 0}$ is a martingale
\cite[Corollary~3.5.13]{KaratzasShreve1991}.

For each fixed $t>0$, Girsanov's theorem
\cite[Theorem~3.5.1]{KaratzasShreve1991}
shows that, under the measure with density $M_t$
on the sigma-algebra $\mathcal F_t$, the process $\{Y_s\}_{0\leqslant s\leqslant t}$
has drift $-x+\theta(x)=-x-R_\k(x)$.
By uniqueness in law, it therefore has the same law
as $\{X_s^\k\}_{0\leqslant s\leqslant t}$ started at $x$.
Consequently, for every bounded measurable function $f$,
\[
\mathbb E_x[f(X_t^\k)]
=
\frac{e^{-\ell t}}{\rho(x)}
\mathbb E_x^{\mathrm{OU}}
\left[
\exp\left(\int_0^tV(Y_s)\,\d s\right)
\rho(Y_t)f(Y_t)
\right].
\]

Disintegrating with respect to $Y_t$ into OU bridge
measures \cite[Proposition~1, Eq.~(2.8)]
{FitzsimmonsPitmanYor1993}, we obtain
\[
p_t^\k(x,y)
=
p_t^{\mathrm{OU}}(x,y)
\frac{\rho(y)}{\rho(x)}e^{-\ell t}
\mathbb E_{x\to y}^{\mathrm{OU},t}
\left[
\exp\left(\int_0^tV(Y_s)\,\d s\right)
\right]
\]
for almost every $y$. Both sides are continuous in $y$.
Indeed, an OU bridge from $x$ to $y$ over $[0,t]$
can be represented as
$Z_s+[x\sinh(t-s)+y\sinh s]/\sinh t$,
where $Z$ is a centered OU bridge from $0$ to $0$
over $[0,t]$ (see \cite{BarczyKern2013}).
Since $V$ is bounded and continuous, dominated
convergence shows that the bridge expectation
is continuous in $y$.
Thus the identity holds for every $y\in\r$.

Solving for the bridge expectation gives
\[
\mathbb E_{x\to y}^{\mathrm{OU},t}
\left[
\exp\left(\int_0^tV(Y_s)\,\d s\right)
\right]
=
e^{\ell t}\frac{\rho(x)}{\rho(y)}
\frac{p_t^\k(x,y)}{p_t^{\mathrm{OU}}(x,y)}.
\]
Substituting \eqref{eqn:p_t^X} and
$\rho=H_{\k-1}/H_\k$ proves \eqref{eq:OU-bridge-functional}.
\qed

\vspace{0.25cm}
\noindent
\textbf{Proof of Proposition~\ref{prop:occupation-time}.}
Write $\ell=2m$, and let $e_1<\cdots<e_m$ and
$o_1<\cdots<o_m$ be the even and odd indices in $\k$,
respectively. We will evaluate $F_{\k,\lambda}(0)$
and $F_{\k,\lambda}'(0)$ directly.

Suppose that $g$ satisfies $g''-2xg'=2\alpha g$.
Differentiating this equation $r$ times and setting $x=0$ gives
$g^{(r+2)}(0)=2(\alpha+r)g^{(r)}(0)$.
Consequently,
\[
g^{(2r)}(0)=p_r(\alpha)g(0),
\qquad
g^{(2r+1)}(0)=q_r(\alpha)g'(0),
\]
where
\[
p_r(z):=2^r\prod_{s=0}^{r-1}(z+2s),
\qquad
q_r(z):=2^r\prod_{s=0}^{r-1}(z+2s+1).
\]
Both $p_r$ and $q_r$ are polynomials of degree $r$
with leading coefficient $2^r$.
These identities apply to $g=F_\lambda$ with
$\alpha=\lambda$ and to $g=H_{k_j}$ with
$\alpha=-k_j$.

Consider the following $(m+1)\times(m+1)$ matrices:
\[
M_{e,F}:=
\begin{pmatrix}
H_{e_1}(0)&\cdots&H_{e_m}(0)&F_\lambda(0)\\
H_{e_1}''(0)&\cdots&H_{e_m}''(0)&F_\lambda''(0)\\
\vdots&&\vdots&\vdots\\
H_{e_1}^{(2m)}(0)&\cdots&H_{e_m}^{(2m)}(0)&F_\lambda^{(2m)}(0)
\end{pmatrix}
\]
and
\[
M_{o,F}:=
\begin{pmatrix}
H_{o_1}'(0)&\cdots&H_{o_m}'(0)&F_\lambda'(0)\\
H_{o_1}^{(3)}(0)&\cdots&H_{o_m}^{(3)}(0)&F_\lambda^{(3)}(0)\\
\vdots&&\vdots&\vdots\\
H_{o_1}^{(2m+1)}(0)&\cdots&H_{o_m}^{(2m+1)}(0)&F_\lambda^{(2m+1)}(0)
\end{pmatrix}.
\]
Let $M_e$ and $M_o$ denote their upper-left
$m\times m$ submatrices, respectively.

We now explain how these matrices arise.
First we observe that $H_{e_j}^{(r)}(0)=0$ for odd $r$, while 
$H_{o_j}^{(r)}(0)=0$ for even $r$. 
Rearranging the rows and columns of the Wronskian defining $H_\k(0)$ gives a block
diagonal matrix of the form
\[
H_{\k}(0)=\epsilon_1 \times  \begin{vmatrix}
M_e & 0\\
0 & M_o
\end{vmatrix},
\]
for some $\epsilon_1\in \{-1,1\}$. 

Similarly, rearranging the rows and columns
of the Wronskian defining $F_{\k,\lambda}(0)$
gives a block lower triangular matrix of the form
\[
F_{\k,\lambda}(0)=\epsilon_2 \times \begin{vmatrix}
M_{e,F} & 0\\
\left[\,0_{m\times m}\ \mathbf v\,\right] & M_o
\end{vmatrix},
\qquad
\mathbf v=
\begin{pmatrix}
F_\lambda'(0)\\
F_\lambda^{(3)}(0)\\
\vdots\\
F_\lambda^{(2m-1)}(0)
\end{pmatrix},
\]
for some $\epsilon_2\in \{-1,1\}$. 
Importantly, we have $\epsilon_1=\epsilon_2$. Indeed, compared with the first
determinant $H_{\k}(0)$, transforming the second determinant $F_{\k,\lambda}(0)$ requires moving the
additional row of even derivatives past $m$
odd-derivative rows, and moving the column
corresponding to $F_\lambda$ past $m$
odd-Hermite columns. Each permutation contributes
a factor $(-1)^m$, so the two factors cancel.

We conclude from the above two formulas that
\begin{equation}\label{eq:Wr-block-ratios1}
\frac{F_{\k,\lambda}(0)}{H_\k(0)}
=\frac{\det M_{e,F}}{\det M_e}.
\end{equation}

To compute $F_{\k,\lambda}'(0)$, we first observe that the
derivative of the Wronskian in \eqref{def_F_k_lambda} is obtained by differentiating only its last row. Indeed, the derivative of a determinant is the sum of the determinants obtained by differentiating one row
at a time and leaving all other rows unchanged. Differentiating any row of a Wronskian determinant other than the last makes it identical to the next row, so the corresponding determinant vanishes.
Thus $F_{\k,\lambda}'(0)$ is represented by
the same determinant as $F_{\k,\lambda}(0)$,
except that its last derivative row has order
$2m+1$ instead of $2m$. Rearranging the rows and columns gives a block
upper triangular matrix of the form
\[
F_{\k,\lambda}'(0)=\epsilon_3\times
\begin{vmatrix}
M_e & \left[\,0_{m\times m}\ \mathbf u\,\right]\\
0 & M_{o,F}
\end{vmatrix},
\qquad
\mathbf u=
\begin{pmatrix}
F_\lambda(0)\\
F_\lambda''(0)\\
\vdots\\
F_\lambda^{(2m-2)}(0)
\end{pmatrix},
\]
for some $\epsilon_3\in\{-1,1\}$. In this case, we also have $\epsilon_3=\epsilon_1$,
since the additional row of odd derivatives
and the column corresponding to $F_\lambda$
are both placed at the end of their respective
odd-derivative and odd-Hermite blocks.
Thus no additional transpositions are needed
compared with the determinant defining $H_\k(0)$.

We conclude that
\begin{equation}\label{eq:Wr-block-ratios2}
\frac{F_{\k,\lambda}'(0)}{H_\k(0)}
=\frac{\det M_{o,F}}{\det M_o}.
\end{equation}

It remains to evaluate the determinant ratios in
\eqref{eq:Wr-block-ratios1} and
\eqref{eq:Wr-block-ratios2}.
Recall that $p_r$ and $q_r$ are polynomials
of degree $r$ with leading coefficient $2^r$.
Therefore, the Vandermonde determinant formula gives,
for any $z_1,\ldots,z_s$,
\begin{equation}\label{eq:pq-Vandermonde}
\det[p_{r-1}(z_j)]_{1\leqslant r,j\leqslant s}
=\det[q_{r-1}(z_j)]_{1\leqslant r,j\leqslant s}
=2^{s(s-1)/2}\prod_{1\leqslant i<j\leqslant s}(z_j-z_i).
\end{equation}
Indeed, elementary row operations reduce either
determinant to the usual Vandermonde determinant
with row $r$ multiplied by $2^{r-1}$.

We first evaluate $\det M_{e,F}/\det M_e$.
Using the identities for the even derivatives at zero,
we factor $H_{e_j}(0)$ from the $j$th column of
$M_{e,F}$, for $1\leqslant j\leqslant m$, and $F_\lambda(0)$
from its last column. The remaining determinant
has entries $p_r(z_j)$, where
\[
(z_1,\ldots,z_{m+1})=(-e_1,\ldots,-e_m,\lambda).
\]
Similarly, factoring $H_{e_j}(0)$ from the
$j$th column of $M_e$ leaves the determinant
with entries $p_r(-e_j)$.

Applying \eqref{eq:pq-Vandermonde} to these
two determinants, we obtain
\[
\det M_{e,F}
=2^{m(m+1)/2}F_\lambda(0)
\Big[\prod_{j=1}^m H_{e_j}(0)\Big]
\Big[\prod_{1\leqslant i<j\leqslant m}(e_i-e_j)\Big]
\prod_{j=1}^m(\lambda+e_j)
\]
and
\[
\det M_e
=2^{m(m-1)/2}
\Big[\prod_{j=1}^m H_{e_j}(0)\Big]
\prod_{1\leqslant i<j\leqslant m}(e_i-e_j).
\]
Since the indices $e_j$ are distinct and
$H_{e_j}(0)\ne0$, we have $\det M_e\ne0$.
Dividing the above identities gives
\begin{equation}\label{eq:ratio-Me}
\frac{\det M_{e,F}}{\det M_e}
=2^mF_\lambda(0)\prod_{j=1}^m(\lambda+e_j).
\end{equation}

The calculation for the odd indices is analogous.
We factor $H_{o_j}'(0)$ from the $j$th column
of $M_{o,F}$ and $F_\lambda'(0)$ from its last
column. The remaining determinant has entries
$q_r(z_j)$, where
\[
(z_1,\ldots,z_{m+1})=(-o_1,\ldots,-o_m,\lambda).
\]
Applying \eqref{eq:pq-Vandermonde}, we obtain
\[
\det M_{o,F}
=2^{m(m+1)/2}F_\lambda'(0)
\Big[\prod_{j=1}^m H_{o_j}'(0)\Big]
\Big[\prod_{1\leqslant i<j\leqslant m}(o_i-o_j)\Big]
\prod_{j=1}^m(\lambda+o_j)
\]
and
\[
\det M_o
=2^{m(m-1)/2}
\Big[\prod_{j=1}^m H_{o_j}'(0)\Big]
\prod_{1\leqslant i<j\leqslant m}(o_i-o_j).
\]
Since the indices $o_j$ are distinct and
$H_{o_j}'(0)\ne0$, we have $\det M_o\ne0$.
Consequently,
\begin{equation}\label{eq:ratio-Mo}
\frac{\det M_{o,F}}{\det M_o}
=2^mF_\lambda'(0)\prod_{j=1}^m(\lambda+o_j).
\end{equation}

Combining \eqref{eq:Wr-block-ratios1},
\eqref{eq:Wr-block-ratios2},
\eqref{eq:ratio-Me} and \eqref{eq:ratio-Mo},
we conclude that
\[
\frac{F_{\k,\lambda}'(0)}{F_{\k,\lambda}(0)}
=\frac{F_\lambda'(0)}{F_\lambda(0)}
\frac{\prod_{j=1}^m(\lambda+o_j)}
{\prod_{j=1}^m(\lambda+e_j)}.
\]
All denominators are nonzero for $\lambda>0$,
since $F_\lambda(0)>0$ and $e_j,o_j>0$.

Finally, recall that
$f_\lambda^+=F_{\k,\lambda}/H_{\k-1}$.
Since $H_{\k-1}$ is even and strictly positive,
we have $H_{\k-1}'(0)=0$, and hence
\[
\frac{(f_\lambda^+)'(0)}{f_\lambda^+(0)}
=\frac{F_{\k,\lambda}'(0)}{F_{\k,\lambda}(0)}.
\]
The standard formulas for parabolic cylinder functions
\cite[Eqs.~(12.2.5)--(12.2.7)]{NIST} give
\[
F_\lambda(0)
=\frac{2^{-\lambda/2}\sqrt{\pi}}
{\Gamma((\lambda+1)/2)},
\qquad
F_\lambda'(0)
=\frac{2^{1-\lambda/2}\sqrt{\pi}}
{\Gamma(\lambda/2)}.
\]
Therefore,
\[
\frac{(f_\lambda^+)'(0)}{f_\lambda^+(0)}
=2\frac{\Gamma((\lambda+1)/2)}{\Gamma(\lambda/2)}
\frac{\prod_{j:\,k_j\ \mathrm{odd}}(\lambda+k_j)}
{\prod_{j:\,k_j\ \mathrm{even}}(\lambda+k_j)}
=2\Phi_\k(\lambda),
\]
which proves the proposition.
\qed

\section*{AI Disclosure}

During the development of this work, we used ChatGPT (OpenAI) to explore mathematical ideas, assist in deriving results, and develop preliminary versions of some proofs. All arguments and proofs presented in the final manuscript were subsequently written and carefully checked by the authors. We take full responsibility for the correctness of the mathematical results and arguments.

%***************************************************************
%***************************************************************
%***************************************************************

%\bibliographystyle{abbrv}
%\bibliography{references}
%***************************************************************
%***************************************************************
%***************************************************************
\end{document}